\documentclass[preprint,12pt]{elsarticle}

\usepackage{graphicx,color}

\usepackage{amssymb}
\usepackage{amsthm}

\usepackage{amsmath,mathrsfs,bm,subfigure,slashbox}

\newtheorem{thm}{Theorem}
\newtheorem{lem}[thm]{Lemma}
\newtheorem{prop}[thm]{Proposition}
\newtheorem{cor}[thm]{Corollary}
\newdefinition{defn}{Definition}
\newdefinition{rmk}{Remark}
\newdefinition{alg}{Algorithm}
\newdefinition{exmp}{Example}
\newproof{pf}{Proof}

\begin{document}

\begin{frontmatter}




\title{Bivariate Quasi-Tower Sets and Their Associated Lagrange Interpolation Bases}
\tnotetext[fund]{This work was supported in part by the National
Grand Fundamental Research 973 Program of China (No. 2004CB318000).}


\author{Tian Dong, Xiaoying Wang}
\author{Shugong Zhang\corref{cor1}}
\cortext[cor1]{Corresponding author} \ead{sgzh@jlu.edu.cn}
\author{Peng Li}
\address{School of Mathemathics,
Key Lab. of Symbolic Computation and Knowledge Engineering
\textup{(}Ministry of Education\textup{)}, Jilin University,
Changchun 130012, PR China}

\begin{abstract}

As we all known, there is still a long way for us to solve arbitrary
multivariate Lagrange interpolation in theory. Nevertheless, it is
well accepted that theories about Lagrange interpolation on special
point sets should cast important lights on the general solution. In
this paper, we propose a new type of bivariate point sets,
quasi-tower sets, whose geometry is more natural than some known
point sets such as cartesian sets and tower sets. For bivariate
Lagrange interpolation on quasi-tower sets, we construct the
associated degree reducing interpolation monomial and Newton bases
w.r.t. common monomial orderings theoretically. Moreover, by
inputting these bases into Buchberger-M\"{o}ller algorithm, we
obtain the reduced Gr\"{o}bner bases for vanishing ideals of
quasi-tower sets much more efficiently than before.

\end{abstract}


\end{frontmatter}


\section{Introduction}\label{sec:int}

Given a set $\Xi=\{\xi^{(1)}, \ldots, \xi^{(\mu)}\} \subset
\mathbb{F}^2$ of $\mu$ distinct points with $\mathbb{F}$ a field.
For prescribed values $f_i \in \mathbb{F}, i=1, \ldots, \mu$,
\emph{bivariate Lagrange interpolation} is to find all polynomials
$p \in \Pi^2$ satisfying
\begin{equation}\label{e:inter}
p(\xi^{(i)}) = f_i,\quad i=1,\ldots,\mu,
\end{equation}
where $\Pi^2:=\mathbb{F}[x, y]$ is the bivariate polynomial ring
over $\mathbb{F}$.

In most cases, for fixed monomial ordering $\prec$, only the bases
for the \emph{degree reducing interpolation space}
$\mathcal{P}\subset \Pi^2$(see \cite{Sau2006}) are pursued for
Lagrange interpolation (\ref{e:inter}). Let $\{p_1, \ldots, p_\mu\}$
be a basis for $\mathcal{P}$ satisfying $p_1 \prec p_2\prec \cdots
\prec p_\mu$. If
$$
p_j(\xi^{(i)})=\delta_{ij},\quad 1\leq i \leq j \leq \mu,
$$
for certain reordering of $\Xi$, then we call $\{p_1, \ldots,
p_\mu\}$ a \emph{degree reducing interpolation Newton basis}(DRINB)
w.r.t. $\prec$ for (\ref{e:inter}). Moreover,
$\mathcal{N}_\prec(\Xi)$, the \emph{Gr\"{o}bner \'{e}scalier} of
$\mathcal{I}(\Xi)$ w.r.t. $\prec$(see \cite{Mor2009}), is the
\emph{degree reducing interpolation monomial basis}(DRIMB) w.r.t.
$\prec$ for (\ref{e:inter}).

Let $\mathcal{I}(\Xi)\subset \Pi^2$ be the vanishing ideal of $\Xi$
and $\mathcal{G}_\prec(\Xi)$ the reduced Gr\"{o}bner basis for
$\mathcal{I}(\Xi)$ w.r.t.$\prec$. In 1982, \cite{MB1982} proposed
Buchberger-M\"{o}ller(BM for short) algorithm for computing bases of
vanishing ideals. Specifically, with inputted $\Xi$ and $\prec$, BM
algorithm outputs $\mathcal{G}_\prec(\Xi)$,
$\mathcal{N}_\prec(\Xi)$, and a DRINB for Lagrange interpolation
(\ref{e:inter}), namely arbitrary bivariate Lagrange interpolation
can be solved with it algorithmly. Nonetheless, the poor complexity
of BM algorithm(see \cite{Lun2008}) severely constrains its
applications.

Since the geometry of $\Xi$ is a dominant factor for the structures
of $\mathcal{G}_\prec(\Xi)$ and $\mathcal{N}_\prec(\Xi)$, theories
about Lagrange interpolation on point sets with describable special
geometries will certainly throw light on the general theoretical
solution of (\ref{e:inter}). \cite{Sau2004, CDZ2006} studied
multivariate Lagrange interpolation on cartesian sets(aka lower
point sets) and \cite{DCZW2009} investigated bivariate Lagrange
interpolation on tower sets.

This paper will introduce, in Section \ref{sec:quasi}, a new type of
bivariate point sets, quasi-tower sets, that is more natural than
cartesian and tower sets. As the section of main results, Section
\ref{sec:main} will discover the degree reducing interpolation bases
w.r.t. certain common monomial ordering such as graded lex order for
bivariate Lagrange interpolation on quasi-tower sets theoretically.
Finally, Section \ref{sec:application} will give an algorithm for
computing reduced Gr\"{o}bner bases for vanishing ideals of
quasi-tower sets.

The following section will serve as a preparation for the paper.

\section{Preliminary}\label{sec:pre}

Let $\mathbb{N}_0$ and $\mathbb{T}^2$ be the monoid of nonnegative
integers and bivariate monomials in $\Pi^2$ respectively. A monomial
in $\mathbb{T}^2$ has the form
$X^{\bm{\alpha}}:=x^{\alpha_1}y^{\alpha_2}$ for
$\bm{\alpha}=(\alpha_1,\alpha_2) \in \mathbb{N}_0^2$.

As in \cite{CLO2007}, we use $\prec_\text{lex}, \prec_\text{invlex},
\prec_\text{grlex}$, and $\prec_\text{grevlex}$ to represent
lexicographic order, inverse lexicographic order, graded lex order,
and graded reverse lex order respectively. Hereafter, we always
assume that $y \prec x$ for any fixed monomial ordering $\prec$.

For a nonzero $f \in \Pi^2$, we let $\mathrm{LM}_\prec(f)$ and
$\bm{\delta_\prec(f)}$ signify the \emph{leading monomial} and
\emph{leading bidegree}($\mathrm{LM}(f)=X^{\bm{\delta_\prec(f)}}$)
of $f$ w.r.t. $\prec$ respectively. Furthermore, set
$$
f\prec g:=\bm{\delta_\prec(f)} \prec \bm{\delta_\prec(g)}.
$$

Let $\mathcal{A}$ be an arbitrary finite subset of $\mathbb{N}_0^2$.
We say that it is \emph{lower} if
$$
\mathrm{R}(\bm{\alpha}):= \{(\alpha'_1,\alpha'_2)\in \mathbb{N}_0^2:
0\leq \alpha'_i \leq \alpha_i, i=1,2\} \subset \mathcal{A}
$$
always holds for any $\bm{\alpha}=(\alpha_1,\alpha_2) \in
\mathcal{A}$. When $\mathcal{A}$ is lower, we let
$m_j=\max_{(h,j)\in \mathcal{A}}h, 0\leq j\leq \nu$, with
$\nu=\max_{(0,k)\in \mathcal{A}} k$. Since $\mathcal{A}$ can be
determined uniquely by $(\nu+1)$-tuple $(m_0,m_1, \ldots,m_{\nu})$,
we represent it as
$$
\mathcal{A}=\mathrm{L}_x(m_0,m_1,\ldots,m_{\nu}).
$$
Likewise, it can also be represented as $\mathrm{L}_y(n_0, \ldots,
n_{m_0})$ with $n_i=\max_{(i,k)\in \mathcal{A}}k$, $0\leq i \leq
m_0$.

As in \cite{DCZW2009}, we constructed two particular lower sets
$S_x(\Xi), S_y(\Xi)\subset \mathbb{N}_0^2$ from $\Xi$. In more
details, we cover $\Xi$ by lines $l_0^x, l_1^x,\ldots, l_\nu^x$
parallel to the $x$-axis and assume that, without loss of
generality, there are $m_j+1$ points, say $u_{0j}^x,
u_{1j}^x,\ldots, u_{m_j,j}^x$, on $l_j^x$ with $m_0\geq m_1\geq
\cdots\geq m_\nu\geq 0$. Set
$$
  S_x(\Xi):=\{(i,j) : 0\leq i \leq m_j,\ 0\leq j\leq \nu\},
$$
which equals to $\mathrm{L}_x(m_0,m_1,\ldots,m_\nu)$. In like
fashion, we cover $\Xi$ by lines $l_0^y, l_1^y,\ldots, l_\lambda^y$
parallel to the $y$-axis and denote the points on $l_i^y$ by
$u_{i0}^y, \ldots, u_{i,n_i}^y$ with $n_0\geq \cdots\geq
n_\lambda\geq 0$. We have
$$
S_y(\Xi):=\{(i,j) : 0\leq i \leq \lambda,\ 0\leq j\leq
n_i\}=\mathrm{L}_y(n_0, n_1, \ldots, n_\lambda).
$$
In addition, the sets of abscissae and ordinates are defined as
$$
    \begin{aligned}
H_j(\Xi):=&\{\bar{x}: (\bar{x}, \bar{y}) \in l_j^x\cap \Xi\},\quad 0\leq j \leq \nu,\\
V_i(\Xi):=&\{\bar{y}: (\bar{x}, \bar{y}) \in l_i^y\cap \Xi\},\quad
0\leq i \leq \lambda.
    \end{aligned}
$$

A point set $\Xi$ satisfying $S_x(\Xi)=S_y(\Xi)$ is called a
cartesian set that has the following property:

\begin{prop}\textup{\cite{DCZW2009}}\label{p:carte}
A point set $\Xi\subset \mathbb{F}^2$ is cartesian if and only if
$$
  H_0(\Xi) \supseteq H_1(\Xi) \supseteq \cdots \supseteq H_\nu(\Xi)
$$
or
$$
  V_0(\Xi) \supseteq V_1(\Xi) \supseteq \cdots \supseteq V_\lambda(\Xi).
$$
\end{prop}

Enlightened by the notion of cartesian sets, \cite{DCZW2009}
introduced tower sets in $\mathbb{F}^2$ as follows.

\begin{defn}\cite{DCZW2009}\label{d:tower}
  We say that a set $\Xi$ of distinct points in $\mathbb{F}^2$ is $x$-\emph{tower} if lower set
$S_x(\Xi)=\mathrm{L}_x(m_0, m_1, \ldots, m_\nu)\subset
\mathbb{N}_0^2 $ with $m_0 > m_1
> \cdots > m_\nu\geq 0$ such that
$$
    \Xi := \{(x_{ij},y_j): (i,j)\in S_x(\Xi) \},
$$
where $x_{ij}\in H_0(\Xi), (i,j)\in S_x(\Xi)$, are distinct for
fixed $j$. Its full name is $S_x(\Xi)$-$x$-\emph{tower} set.
Similarly, if lower set $S_y(\Xi)=\mathrm{L}_y(n_0, n_1, \ldots,
n_\lambda)\subset \mathbb{N}_0^2$, $n_0 > n_1
> \cdots > n_\lambda\geq 0$, such that
$$
    \Xi := \{(x_i,y_{ij}): (i,j)\in S_y(\Xi)\},
$$
where $y_{ij}\in V_0(\Xi), (i,j)\in S_y(\Xi)$, are distinct for
fixed $i$, we will call $\Xi$ a $y$-\emph{tower} set(or
$S_y(\Xi)$-$y$-\emph{tower} set in its full name).
\end{defn}

About the bivariate Lagrange interpolation on a tower set,
\cite{DCZW2009} proved the succeeding theorem.

\begin{thm}\textup{\cite{DCZW2009}}\label{t:td}
Given an $x$-tower set $\Xi\subset \mathbb{F}^2$. The DRIMB w.r.t.
$\prec_{\mathrm{grlex}}$ or $\prec_{\mathrm{lex}}$ for
\eqref{e:inter} is
$$
\{x^iy^j: (i, j) \in S_x(\Xi)\}.
$$
If $\Xi$ is $y$-tower, then the DRIMB w.r.t.
$\prec_{\mathrm{grevlex}}$ or $\prec_{\mathrm{invlex}}$ for
\eqref{e:inter} is
$$
\{x^iy^j: (i, j) \in S_y(\Xi)\}.
$$
\end{thm}

\section{Quasi-Tower Sets}\label{sec:quasi}

In this section, we will introduce the notion of generalized tower
and quasi-tower sets and compare them with tower sets and cartesian
sets.

\begin{defn}\label{d:gtower}
  Let $\Xi$ be a set of distinct points in $\mathbb{F}^2$. We call
  $\Xi$ a \emph{generalized $x$-tower} set if $S_x(\Xi)=\mathrm{L}_x(m_0, m_1, \ldots, m_\nu)\subset
  \mathbb{N}_0^2 $ with $m_0 > m_1 > \cdots > m_\nu\geq 0$ such that
$$
    \Xi := \{(x_{ij},y_j): (i,j)\in S_x(\Xi) \},
$$
where $x_{ij}$ are distinct for fixed $j$. In like fashion, if lower
set $S_y(\Xi)=\mathrm{L}_y(n_0, n_1, \ldots, n_\lambda)\subset
\mathbb{N}_0^2$, $n_0 > n_1 > \cdots > n_\lambda\geq 0$, such that
$$
    \Xi := \{(x_i,y_{ij}): (i,j)\in S_y(\Xi)\},
$$
where $y_{ij}$ are distinct for fixed $i$, we will call $\Xi$ a
\emph{generalized $y$-tower} set.

\end{defn}

Observe Definition \ref{d:tower} and \ref{d:gtower}. We find that
they are similar to each other except two more conditions in
Definition \ref{d:tower} which implies that a tower set is always
generalized tower, namely the notion of generalized tower sets is
actually a generalization of tower sets'. Therefore, what interests
us most now should be the generalized tower sets that are not tower.

\begin{defn}\label{d:qtower}
  If $\Xi\subset \mathbb{F}^2$ is a generalized $x$-tower set that is not $x$-tower, then we
  call it a \emph{quasi-$x$-tower} set. Similarly, a non-$y$-tower
  generalized $y$-tower set is called a \emph{quasi-$y$-tower} set.
\end{defn}

From Definition \ref{d:gtower} and \ref{d:qtower} we know that the
geometry of quasi-tower sets is much freer than tower sets' in the
sense that they need not fulfill the two extra conditions that
demand $x_{ij}\subset H_0(\Xi)$ for $x$-tower and $x_{ij}\subset
V_0(\Xi)$ for $y$-tower. The following propositions will show us how
did quasi-tower sets acquire their name.

\begin{prop}\label{p:qx}
  Let $\Xi$ be a quasi-$x$-tower set with $S_x(\Xi)=\mathrm{L}_x(m_0, \ldots, m_\nu)$.
  If $B\subset \mathbb{F}^2$ is a set of distinct horizontal points such that $B\cap \Xi=\emptyset$ and
  $\bigcup_{i=0}^\nu H_i(\Xi)\subsetneq H_0(B)$, then $\Xi\cup B$ is
  $x$-tower and will be known as a \emph{derived $x$-tower} set from $\Xi$ with \emph{base set} $B$.
\end{prop}

\begin{pf}
  By hypothesis, we assume that $S_y(\Xi)=\mathrm{L}_y(n_0, n_1, \ldots, n_\lambda)$ and
  $B=\{(x_0^B,y^B), (x_1^B, y^B), \ldots, (x_{m_B}^B,y^B)\}$ where
  $x_i^B\neq x_j^B, i\neq j$, and
  $y^B\notin \bigcup_{j=0}^\lambda V_j(\Xi)$. Recall the construction process of $S_x(\Xi)$ in
  Section \ref{sec:pre}. Since $\bigcup_{i=0}^\nu H_i(\Xi)\subsetneq H_0(B)$ and $\Xi$ is
  quasi-$x$-tower, we can deduce easily that
   $$
   S_x(\Xi\cup B)=\mathrm{L}_x(m_B, m_0, \ldots, m_\nu),
   $$
   where $m_B>m_0>\cdots>m_\nu\geq 0$. From Definition \ref{d:gtower} and \ref{d:qtower}, $\Xi$ is quasi-$x$-tower
   implies that
   $$
   \Xi=\{(x_{ij}^\Xi, y_j^\Xi): (i, j)\in S_x(\Xi)\},
   $$
   where $x_{ij}^\Xi$ are distinct for fixed $j$, hence, we have
   $$
   \Xi\cup B=\{(x_{ij}, y_j): (i, j)\in S_x(\Xi\cup B)\}
   $$
   where
    \begin{align*}
    x_{i0}&=x_i^B, y_0=y^B, & i&=0, 1, \ldots, m_B,\\
    x_{i,j+1}&=x_{ij}^\Xi, y_{j+1}=y_j^\Xi, & i&=0, 1, \ldots, m_j^\Xi, j=0, 1, \ldots,
    \nu,
   \end{align*}
   and $x_{ij}$ are distinct for fixed $j$. Therefore, by
   Definition \ref{d:tower}, $\Xi\cup B$ is $x$-tower.
   \qed
\end{pf}

In the same way, we can prove the following proposition for
quasi-$y$-tower sets.

\begin{prop}\label{p:qy}
  Let $\Xi$ be a quasi-$y$-tower set with $S_y(\Xi)=\mathrm{L}_y(n_0, n_1, \ldots, n_\lambda)$.
  If $B\subset \mathbb{F}^2$ is a set of distinct vertical points such that $B\cap \Xi=\emptyset$ and
  $\bigcup_{j=0}^\lambda V_j(\Xi)\subsetneq V_0(B)$, then $\Xi\cup B$ is
  $y$-tower and will be named a \emph{derived $y$-tower} set from $\Xi$ with \emph{base set} $B$.
\end{prop}

Proposition \ref{p:qx} and \ref{p:qy} imply the following corollary
immediately.

\begin{cor}\label{c:qexist}
  Let $\Xi\subset \mathbb{F}^2$ be a quasi-tower set. Then there are
  infinite derived tower sets from $\Xi$.
\end{cor}

\begin{figure}[!htbp]
\centering
\subfigure[$\Xi$]{\includegraphics[width=6cm,height=5.3cm]{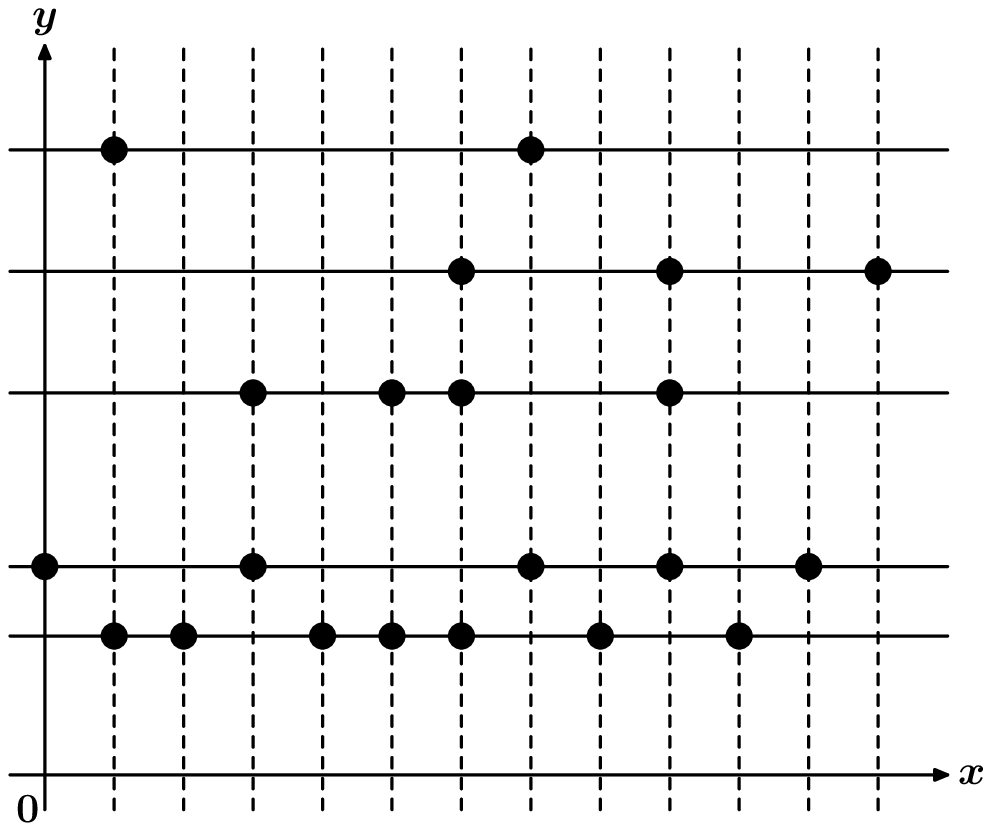}}
\subfigure[$\Xi\cup
B$]{\includegraphics[width=6cm,height=5.3cm]{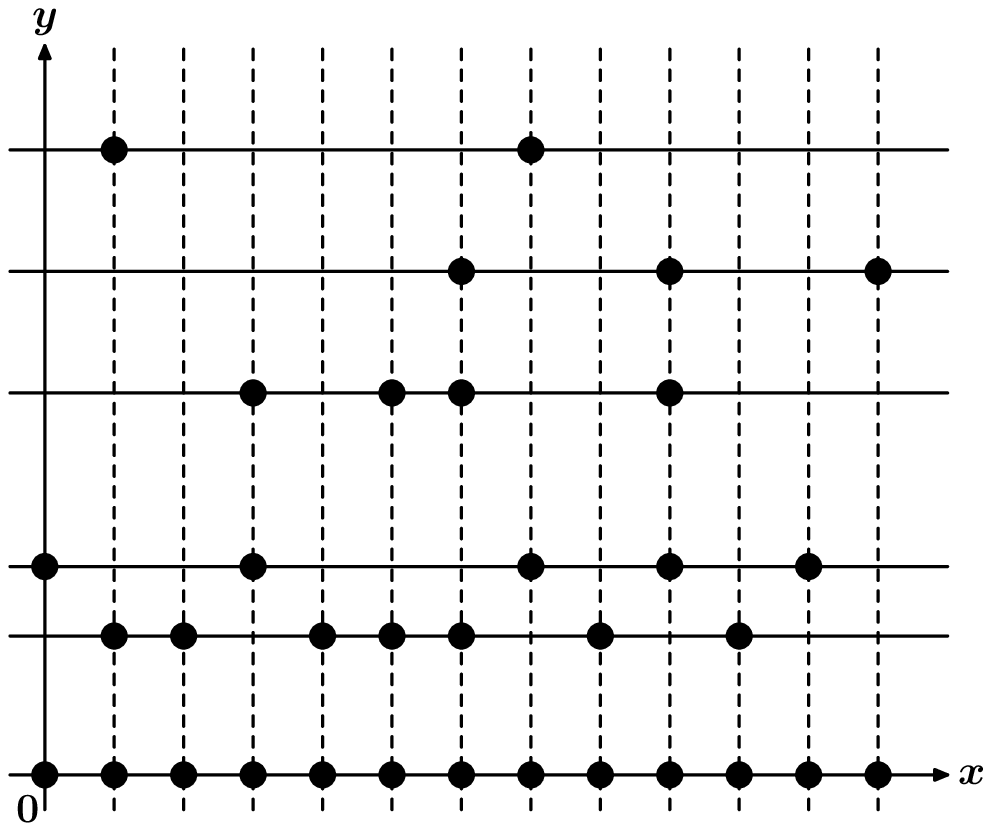}} \caption{A
quasi-$x$-tower set $\Xi$ and a derived $x$-tower set $\Xi\cup B$
from $\Xi$ with base set $B$}\label{f:qtt}
\end{figure}

\begin{exmp}\label{ex:qtt}
\begin{align*}
    \Xi=\{&(0.2,0.4), (0.4,0.4), (0.8,0.4), (1,0.4), (1.2,0.4), (1.6,0.4),\\
          &(2,0.4), (0,0.6), (0.6,0.6), (1.4,0.6), (1.8,0.6), (2.2,0.6), \\
          &(0.6,1.1), (1,1.1),(1.2,1.1), (1.8,1.1), (1.2,1.45), (1.8,1.45), \\
          &(2.4,1.45), (0.2,1.8), (1.4,1.8)\}\subset \mathbb{Q}^2
\end{align*}
is a quasi-$x$-tower set(shown in (a) of Figure \ref{f:qtt}) and
\begin{align*}
    B=\{&(0,0), (0.2,0), (0.4,0), (0.6,0), (0.8,0), (1,0),\\
        &(1.2,0), (1.4,0), (1.6,0), (1.8,0), (2,0), (2.2,0),
        (2.4,0)\}.
\end{align*}
Obviously, $\Xi\cup B$(shown in (b) of Figure \ref{f:qtt}) is an
$x$-tower set that is derived from $\Xi$ with base set $B$.

\end{exmp}

Proposition 6 of \cite{DCZW2009} discovered the relation between
cartesian sets and tower sets. The next proposition shows the
relation between cartesian sets and quasi-tower sets.

\begin{prop}
  A quasi-tower set is not a cartesian set and vice versa.
\end{prop}

\begin{pf}
  From Definition \ref{d:qtower} we can deduce easily that a point
set $\Xi\subset \mathbb{F}^2$ is quasi-$x$-tower if and only if it
is a generalized $x$-tower set such that there exists at least one
point $\xi\in \Xi$ whose abscissa is not in $H_0(\Xi)$. Proposition
\ref{p:carte} implies the statement of this proposition immediately
for quasi-$x$-tower sets. The quasi-$y$-tower cases can be proved
similarly. \qed
\end{pf}

\section{Bases for Lagrange Interpolation on Quasi-Tower
Sets}\label{sec:main}

In this section, we will present the degree reducing interpolation
bases for Lagrange interpolation on quasi-tower sets w.r.t.
$\prec_{\mathrm{lex}}, \prec_{\mathrm{grevlex}}$ etc. theoretically.
Let us begin with the following lemmas.

\begin{lem}\label{l:line}
  Fix monomial ordering $\prec_{\mathrm{grlex}}$. Let
  $$
  \Xi'=\{(x_{ij}, y_j): (i, j)\in S_x(\Xi')\}
  $$
  be an $x$-tower set in $\mathbb{F}^2$ with $S_x(\Xi')=\mathrm{L}_x(m_0, m_1, \ldots,
  m_\nu)$.
  If $\mathcal{G}_{\prec_{\mathrm{grlex}}}(\Xi')=\{g_1, g_2, \ldots, g_s\}$
  with $g_s \prec_{\mathrm{grlex}} \cdots \prec_{\mathrm{grlex}} g_1$, then
  \begin{enumerate}
    \item[\textup{(i)}] $s=\nu+2$.
    \item[\textup{(ii)}] $g_1=(x-x_{0,0})(x-x_{1,0})\cdots(x-x_{m_0,0})$.
    \item[\textup{(iii)}] $g_s=(y-y_0)(y-y_1)\cdots(y-y_\nu)$.
    \item[\textup{(iv)}] $g_j\in\langle y-y_0 \rangle$ and
     $$
     \mathrm{LM}_{\prec_{\mathrm{lex}}}(g_j)=\mathrm{LM}_{\prec_{\mathrm{grlex}}}(g_j)
     =x^{m_{j-1}+1}y^{j-1},\quad j=2,\ldots,s-1.
     $$
    \item[\textup{(v)}] $\mathcal{I}(\{(x_{0,0}, y_0), \ldots, (x_{m_{0},0}, y_0)\})=\langle g_1,y-y_0\rangle.$
  \end{enumerate}
\end{lem}

\begin{pf}
From Theorem \ref{t:td} we know that
$\mathcal{N}_{\prec_{\mathrm{grlex}}}(\Xi')=\{x^iy^j: (i, j)\in
S_x(\Xi')\}$. Since $S_x(\Xi')=\mathrm{L}_x(m_0, m_1, \ldots,
  m_\nu)$, it is easily checked that
  $$
  \{\mathrm{LM}_{\prec_{\mathrm{grlex}}}(g_i): 1\leq i \leq s\}
  =\{x^{m_j+1}y^j: 0\leq j \leq \nu\}\cup\{y^{\nu+1}\},
  $$
 which implies (i) readily. In the proof of Theorem
  \ref{t:td}(For brevity, we refer the reader to \cite{DCZW2009} for
  details), we find that $(x-x_{0,0})(x-x_{1,0})\cdots(x-x_{m_0,0}),
  (y-y_0)(y-y_1)\cdots(y-y_\nu)\in
  \mathcal{G}_{\prec_{\mathrm{grlex}}}(\Xi')$.
  Since $\Xi'$ is
  $x$-tower, we have $m_0>m_1>\cdots>m_\nu\geq 0$. Thus, for $0\leq j\leq \nu-1$,
  \begin{align*}
    \deg x^{m_j+1}y^j&=m_j+1+j\\
    &\geq m_{j+1}+1+j+1\\
    &=\deg x^{m_{j+1}+1}y^{j+1}.
  \end{align*}
  Next, we claim that $\deg x^{m_j+1}y^j \geq \deg y^{\nu+1}$, that is $m_j\geq \nu-j$,
  $j=0,\ldots,\nu$. In fact, it is easy to see from $m_0>m_1>\cdots>m_\nu\geq 0$
  that $m_j-m_\nu\geq \nu-j$, hence, $m_j\geq \nu-j, j=0, \ldots, \nu$.
  Therefore, (ii), (iii) follows.

  Still from the proof of Theorem
  \ref{t:td}, we know that $g_j\in\langle y-y_0 \rangle$ and there is no monomial
   $x^{\alpha_1} y^{\alpha_2}$ of $g_j$ except $\mathrm{LM}_{\prec_{\mathrm{grlex}}}(g_j)$
   satisfying $\alpha_1\geq m_{j-1}+1$. Therefore, the
   definition of $\prec_{\mathrm{lex}}$ implies (iv).

Finally, (v) is a direct consequence of Lemma 8 in \cite{DCZW2009}.
\qed
\end{pf}

\begin{lem}\textup{\cite{CLO2007}}\label{l:idqu}
Let $\Xi$ and $\Xi'$ be varieties in $\mathbb{F}^d$. Then
$$
    \mathcal{I}(\Xi):\mathcal{I}(\Xi')=\mathcal{I}(\Xi-\Xi').
$$

\end{lem}

%

\begin{lem}\textup{\cite{CLO2007}}\label{l:idqu3}
Let $\mathcal{I}$ be an ideal and $g$ an element of $\Pi^2$. If
$\{h_1, \ldots, h_p\}$ is a basis of the ideal $\mathcal{I}\cap
\langle g\rangle$, then $\{h_1/g, \ldots, h_p/g\}$ is a basis of
$\mathcal{I}: \langle g\rangle$.

\end{lem}

The following theorem shows us the degree reducing interpolation
bases for bivariate Lagrange interpolation on quasi-$x$-tower sets
w.r.t. $\prec_{\mathrm{grlex}}$.

\begin{thm}\label{t:qtowerx}
  Let $\Xi$ be a quasi-$x$-tower set in $\mathbb{F}^2$.
  The DRIMB for Lagrange interpolation on $\Xi$ w.r.t. $\prec_{\mathrm{grlex}}$ is
$$
    \mathcal{N}_{\prec_{\mathrm{grlex}}}(\Xi)=\{x^iy^j: (i, j)\in
    S_x(\Xi)\}.
$$
\end{thm}

\begin{pf}
  Firstly, we assume that $S_x(\Xi)=\mathrm{L}_x(m_1, \ldots, m_\nu)$. By Corollary \ref{c:qexist},
  there exists a derived $x$-tower set from $\Xi$, say $\Xi'$, with base set
  $B=\{(x_0,y_0), (x_1, y_0), \ldots, (x_{m_0},y_0)\}$. Thus,
  Proposition \ref{p:qx} implies that $S_x(\Xi')=\mathrm{L}_x(m_0, m_1, \ldots, m_\nu)$. Moreover, it follows
  from (v) of Lemma \ref{l:line} that
  \begin{align*}
    \mathcal{I}(B)&=\langle
    (x-x_{0})(x-x_{1})\cdots(x-x_{m_0}),y-y_0\rangle\\
                  &=\langle
    (x-x_{0})(x-x_{1})\cdots(x-x_{m_0})\rangle+\langle y-y_0\rangle.
  \end{align*}
  Since $B\cap \Xi=\emptyset$, Lemma \ref{l:idqu} leads to $\mathcal{I}(\Xi)=\mathcal{I}(\Xi'-B)
  =\mathcal{I}(\Xi'):\mathcal{I}(B)$. Therefore, we have
  \begin{align*}
    \mathcal{I}(\Xi)&=\mathcal{I}(\Xi'): \big(\langle
    (x-x_{0})(x-x_{1})\cdots(x-x_{m_0})\rangle+\langle
    y-y_0\rangle\big)\\
                    &=\big(\mathcal{I}(\Xi'): \langle
   (x-x_{0})(x-x_{1})\cdots(x-x_{m_0})\rangle\big)\cap
   \big(\mathcal{I}(\Xi'): \langle y-y_0\rangle\big).
  \end{align*}
  To be sure, it follows from (ii) of Lemma \ref{l:line} that $\langle
   (x-x_{0})(x-x_{1})\cdots(x-x_{m_0})\rangle\subset
   \mathcal{I}(\Xi')$, hence,
  $$
  \mathcal{I}(\Xi'): \langle
   (x-x_{0})(x-x_{1})\cdots(x-x_{m_0})\rangle=\Pi^2,
  $$
  which means
  $$
    \mathcal{I}(\Xi)=\mathcal{I}(\Xi'): \langle y-y_0\rangle.
  $$
  Recalling Lemma \ref{l:line}, we suppose that $\mathcal{G}_{\prec_{\mathrm{grlex}}}(\Xi')=
  \{g_1, g_2, \ldots, g_{\nu+2}\}$ with $g_{\nu+2} \prec_{\mathrm{grlex}} \cdots \prec_{\mathrm{grlex}} g_1$.
  Hence, we have
  \begin{align*}
    g_1&=(x-x_{0})(x-x_{1})\cdots(x-x_{m_0}),\\
    g_i&\in\langle y-y_0 \rangle,\quad i=2,\ldots, \nu+2.
  \end{align*}
  Next, for the basis for $\mathcal{I}(\Xi)$, we should
  compute the basis for $\mathcal{I}(\Xi')\cap\langle
  y-y_0\rangle$ so that we can apply Lemma \ref{l:idqu3}. Let $t$ be a new variable. According to
  \cite{CLO2007},
  \begin{align*}
    \mathcal{I}(\Xi')\cap\langle y-y_0\rangle&= \langle t\mathcal{I}(\Xi'),
   (1-t)\langle y-y_0\rangle \rangle\cap \Pi^2\\
   &=\langle tg_1, tg_2, \ldots, tg_{\nu+2}, (1-t)(y-y_0) \rangle \cap
   \Pi^2.
  \end{align*}
  Since $g_i\in\langle y-y_0 \rangle, i=2,\ldots, \nu+2$, $g_i/(y-y_0)\in
  \Pi^2$, hence,
  \begin{align*}
    &tg_i+(g_i/(y-y_0))(1-t)(y-y_0)\\
    =&(g_i/(y-y_0))t(y-y_0)+(g_i/(y-y_0))(1-t)(y-y_0)\\
    =&(g_i/(y-y_0))(y-y_0)\\
    =&g_i,
  \end{align*}
  holds for $i=2,\ldots, \nu+2$. Thus,
  $$
    \mathcal{I}(\Xi')\cap\langle y-y_0\rangle=\langle tg_1, (1-t)(y-y_0), g_2, \ldots, g_{\nu+2} \rangle \cap
   \Pi^2.
  $$
  Since $\mathcal{G}_{\prec_{\mathrm{grlex}}}(\Xi')=\{g_1, g_2, \ldots, g_{\nu+2}\}$, it is easily checked
  from Buchberger's S-pair criterion and (iv) of Lemma \ref{l:line} that
  $$
  \{tg_1, (1-t)(y-y_0), g_2, \ldots, g_{\nu+2}\}
  $$
  is the reduced Gr\"{o}bner basis for $\langle t\mathcal{I}(\Xi'),
   (1-t)(y-y_0)\rangle$ w.r.t. $\prec_{\mathrm{lex}}$ with
    $y \prec_{\mathrm{lex}} x \prec_{\mathrm{lex}} t$. Consequently,
  $$
   \mathcal{I}(\Xi')\cap\langle y-y_0\rangle=\langle g_2, \ldots, g_{\nu+2}
   \rangle,
  $$
  hence,
  $$
    \mathcal{I}(\Xi)=\mathcal{I}(\Xi'):\langle y-y_0\rangle= \langle g_2/(y-y_0), \ldots,
    g_{\nu+2}/(y-y_0)\rangle,
  $$
  follows from Lemma \ref{l:idqu3}. In fact, we can deduce that
  \begin{align*}
    \mathcal{G}_{\prec_{\mathrm{grlex}}}(\Xi)&=\{g_2/(y-y_0), \ldots,
    g_{\nu+2}/(y-y_0)\}\\
    &=:\{g_1', \ldots, g_{\nu+1}'\}.
  \end{align*}
  By Lemma \ref{l:line}, $\Xi'$ is $x$-tower implies that
  $$
  \{\mathrm{LM}_{\prec_{\mathrm{grlex}}}(g_j): 1\leq j\leq \nu+2\}
  =\{x^{m_j+1}y^j: 0\leq j \leq \nu\}\cup\{y^{\nu+1}\},
  $$
  therefore,
  $$
  \{\mathrm{LM}_{\prec_{\mathrm{grlex}}}(g_j'): 1\leq j\leq \nu+1\}
  =\{x^{m_j+1}y^{j-1}: 1\leq j \leq \nu\}\cup\{y^{\nu}\}.
  $$
  As a result, since $S_x(\Xi)=\mathrm{L}_x(m_1, \ldots, m_\nu)$, we have
  $$
  \mathcal{N}_{\prec_{\mathrm{grlex}}}(\Xi)=\{x^iy^j: (i, j)\in S_x(\Xi)\}
  $$
  which finishes the proof. \qed
\end{pf}

Theorem \ref{t:qtowerx} has solved the DRIMB problems w.r.t.
$\prec_{\mathrm{grlex}}$. Next, let us turn to DRINB problems.

\begin{thm}\label{t:qtowerphi^x}
Let $\Xi\subset \mathbb{F}^2$ be a quasi-$x$-tower set of points
$$
u_{mn}^x=(x_{mn}, y_n),\quad (m, n)\in S_x(\Xi),
$$
which give rise to polynomials
$$
\phi_{ij}^x=\varphi_{ij}^x
\prod_{t=0}^{j-1}(y-y_{t})\prod_{s=0}^{i-1}(x-x_{sj}), \quad
(i,j)\in S_x(\Xi),
$$
where
$\varphi_{ij}^x=1/\prod_{t=0}^{j-1}(y_{j}-y_{t})\prod_{s=0}^{i-1}(x_{ij}-x_{sj})\in
\mathbb{F}$, and the empty products are taken as 1. Then
\begin{equation*}\label{e:Qx}
Q_x=\{\phi_{ij}^x: (i, j)\in S_x(\Xi)\}
\end{equation*}
is a \emph{DRINB} w.r.t. $\prec_{\mathrm{grlex}}$ for
\eqref{e:inter} satisfying
$$
\phi_{ij}^x(u_{mn}^x)=\delta_{(i,j), (m,n)}, \quad
(i,j)\succeq_{\mathrm{invlex}} (m,n).
$$
\end{thm}

\begin{pf}
Fix $(i, j)\in S_x(\Xi)$. When $(i, j)=(m, n)$, since $x_{0j}\neq
x_{1j}\neq\cdots\neq x_{ij}$ and $y_0\neq y_1\neq\cdots\neq y_j$,
$$
\phi_{ij}^x(u_{ij}^x)=\varphi_{ij}^x\prod_{t=0}^{j-1}(y_j-y_t)\prod_{s=0}^{i-1}(x_{ij}-x_{sj})\\
                     =\varphi_{ij}^x/\varphi_{ij}^x=1
$$
follows. Otherwise, if $(i, j)\succ_{\mathrm{invlex}} (m, n)$, we
have $j>n$, or $j=n, i>m$. When $j>n$, we have
$$
\phi_{ij}^x(u_{mn}^x)=\varphi_{ij}^x(y_n-y_0)\cdots(y_n-y_n)\cdots(y_n-y_{j-1})\prod_{s=0}^{i-1}(x_{mn}-x_{sj})=0,
$$
and when $j=n, i>m$,
\begin{align*}
\phi_{ij}^x(u_{mn}^x)&=\varphi_{ij}^x\prod_{t=0}^{j-1}(y_n-y_t)(x_{mn}-x_{0j})\cdots(x_{mn}-x_{mj})\cdots(x_{mn}-x_{i-1,j})\\
&=\varphi_{ij}^x\prod_{t=0}^{n-1}(y_n-y_t)(x_{mn}-x_{0n})\cdots(x_{mn}-x_{mn})\cdots(x_{mn}-x_{i-1,n})\\
&=0,
\end{align*}
which leads to
$$
\phi_{ij}^x(u_{mn}^x)=0, \quad (i,j)\succ_{\mathrm{invlex}} (m,n),
$$
namely $Q_x$ is a Newton basis for $\mathrm{Span}_\mathbb{F}Q_x$. By
Theorem \ref{t:qtowerx}, it is easy to see that
$\mathrm{Span}_\mathbb{F}Q_x=\mathrm{Span}_\mathbb{F}{\mathcal{N}_{\prec_{\mathrm{grlex}}}(\Xi)}$.
Therefore, $Q_x$ is a DRINB w.r.t. $\prec_{\mathrm{grlex}}$ for
\eqref{e:inter}. \qed
\end{pf}

Similarly, we can prove the following theorems:

\begin{thm}\label{t:qtowery}
  Let $\Xi$ be a quasi-$y$-tower set in $\mathbb{F}^2$.
  The DRIMB for Lagrange interpolation on $\Xi$ w.r.t. $\prec_{\mathrm{grevlex}}$ is
$$
    \mathcal{N}_{\prec_{\mathrm{grevlex}}}(\Xi)=\{x^iy^j: (i, j)\in
    S_y(\Xi)\}.
$$
\end{thm}

\begin{thm}\label{t:qtowerphi^y}
Let $\Xi\subset \mathbb{F}^2$ be a quasi-$y$-tower set of points
$$
u_{mn}^y=(x_m, y_{mn}), (m, n)\in S_y(\Xi).
$$
We define the polynomials
$$
\phi_{ij}^y=\varphi_{ij}^y\prod_{s=0}^{i-1}(x-x_s)\prod_{t=0}^{j-1}(y-y_{it}),\quad
(i,j)\in S_y(\Xi),
$$
where
$\varphi_{ij}^y=1/\prod_{s=0}^{i-1}(x_i-x_s)\prod_{t=0}^{j-1}(y_{ij}-y_{it})\in
\mathbb{F}$. The empty products are taken as 1. Then,
\begin{equation}\label{e:Qy}
Q_y=\{\phi_{ij}^x: (i, j)\in S_y(\Xi)\}
\end{equation}
is a \emph{DRINB} w.r.t. $\prec_{\mathrm{grevlex}}$ for
\eqref{e:inter} satisfying
$$
\phi_{ij}^y(u_{mn}^y)=\delta_{(i,j),(m,n)}, \quad
(i,j)\succeq_{\mathrm{lex}} (m,n).
$$
\end{thm}

So far, we have discovered the degree reducing interpolation bases
for Lagrange interpolation on quasi-tower sets w.r.t.
$\prec_{\mathrm{grlex}}$ and $\prec_{\mathrm{grevlex}}$. Next, we
should turn to $\prec_{\mathrm{lex}}$ and $\prec_{\mathrm{invlex}}$
cases.

\begin{lem}\textup{\cite{WZD2009}}\label{l:WZD}
Let
\begin{align*}
   \Xi&=\{u_{mn}^x=(x_{mn}^x, y_{mn}^x): (m, n)\in S_x(\Xi)\}\\
      &=\{u_{mn}^y=(x_{mn}^y, y_{mn}^y): (m, n)\in S_y(\Xi)\}
\end{align*}
be a set of distinct points in $\mathbb{F}^2$. Then

 \textup{(i)} the set
$\{x^iy^j : (i,j)\in S_x(\Xi)\}$ is the \emph{DRIMB} as well as
$\{\phi_{ij}^x: (i, j)\in S_x(\Xi)\}$ is a \emph{DRINB}  w.r.t.
$\prec_{\mathrm{lex}}$ for \eqref{e:inter}, where
$$
\phi_{ij}^x=\varphi_{ij}^x
\prod_{t=0}^{j-1}(y-y_{0t}^x)\prod_{s=0}^{i-1}(x-x_{sj}^x), \quad
(i,j)\in S_x(\Xi),
$$
with
$\varphi_{ij}^x=1/\prod_{t=0}^{j-1}(y_{0j}^x-y_{0t}^x)\prod_{s=0}^{i-1}(x_{ij}^x-x_{sj}^x)\in
\mathbb{F}$ and the empty products taken as 1;

\textup{(ii)} the set $\{x^iy^j : (i, j)\in S_y(\Xi)\}$ is the
\emph{DRIMB} as well as $\{\phi_{ij}^y: (i, j)\in S_y(\Xi)\}$ is a
\emph{DRINB} w.r.t. $\prec_{\mathrm{invlex}}$ for
\textup{(\ref{e:inter})}, where
$$
\phi_{ij}^y=\varphi_{ij}^y\prod_{s=0}^{i-1}(x-x_{s0}^y)\prod_{t=0}^{j-1}(y-y_{it}^y),\quad
(i,j)\in S_y(\Xi),
$$
with
$\varphi_{ij}^y=1/\prod_{s=0}^{i-1}(x_{i0}^y-x_{s0}^y)\prod_{t=0}^{j-1}(y_{ij}^y-y_{it}^y)\in
\mathbb{F}$, and the empty products are taken as 1.
\end{lem}

Finally, Theorem \ref{t:qtowerx}-\ref{t:qtowerphi^y} and Lemma
\ref{l:WZD} together give rise to our main theorem.

\begin{thm}\label{t:main}
Let
$$
   \Xi^x=\{u_{ij}^x=(x_{ij}^x, y_{j}^x): (i, j)\in S_x(\Xi^x)\}
$$
be a quasi-$x$-tower set in $\mathbb{F}^2$. Then the set $\{x^iy^j :
(i,j)\in S_x(\Xi^x)\}$ is the \emph{DRIMB} as well as
$\{\phi_{ij}^x: (i, j)\in S_x(\Xi^x)\}$ is a \emph{DRINB} w.r.t.
$\prec_{\mathrm{lex}}$ or $\prec_{\mathrm{grlex}}$ for
\eqref{e:inter}, where
$$
\phi_{ij}^x=\varphi_{ij}^x
\prod_{t=0}^{j-1}(y-y_{t}^x)\prod_{s=0}^{i-1}(x-x_{sj}^x), \quad
(i,j)\in S_x(\Xi^x),
$$
with
$\varphi_{ij}^x=1/\prod_{t=0}^{j-1}(y_{j}^x-y_{t}^x)\prod_{s=0}^{i-1}(x_{ij}^x-x_{sj}^x)\in
\mathbb{F}$ and the empty products taken as 1.

If
$$
   \Xi^y=\{u_{ij}^y=(x_{i}^y, y_{ij}^y): (i, j)\in
   S_y(\Xi^y)\}\subset \mathbb{F}^2
$$
is a quasi-$y$-tower set. Then the set $\{x^iy^j : (i, j)\in
S_y(\Xi^y)\}$ is the \emph{DRIMB} as well as $\{\phi_{ij}^y: (i,
j)\in S_y(\Xi^y)\}$ is a \emph{DRINB} w.r.t.
$\prec_{\mathrm{invlex}}$ or $\prec_{\mathrm{grevlex}}$ for
\textup{(\ref{e:inter})}, where
$$
\phi_{ij}^y=\varphi_{ij}^y\prod_{s=0}^{i-1}(x-x_{s}^y)\prod_{t=0}^{j-1}(y-y_{it}^y),\quad
(i,j)\in S_y(\Xi^y),
$$
with
$\varphi_{ij}^y=1/\prod_{s=0}^{i-1}(x_{i}^y-x_{s}^y)\prod_{t=0}^{j-1}(y_{ij}^y-y_{it}^y)\in
\mathbb{F}$, and the empty products are taken as 1.
\end{thm}

\begin{exmp}\label{ex:qtt2}
  Let us continue with Example \ref{ex:qtt}. Observeing Figure
  \ref{f:qtt}, we have
  $$
  S_x(\Xi)=\mathrm{L}_x(6, 4, 3, 2, 1).
  $$
  Therefore, according to Theorem \ref{t:main},
  \begin{equation}\label{e:exN}
    \begin{aligned}
    \mathcal{N}_{\prec_{\mathrm{grlex}}}(\Xi)=\mathcal{N}_{\prec_{\mathrm{lex}}}(\Xi)=&
   \{x^iy^j: (i, j)\in S_x(\Xi)\}\\
   =&\{1, x, x^2, x^3, x^4, x^5, x^6, y, xy, x^2y, x^3y, x^4y, \\
    & \phantom{\{}y^2, xy^2, x^2y^2, x^3y^2, y^4, xy^4, x^2y^4, y^5, xy^5\}
  \end{aligned}
  \end{equation}
  and a DRINB for interpolation on $\Xi$ w.r.t.
$\prec_{\mathrm{lex}}$ or $\prec_{\mathrm{grlex}}$ is
  \begin{equation}\label{e:exQ}
  \begin{aligned}
    &\Bigg\{1, 5(x-0.2), \frac{25}{6}(x-0.2)(x-0.4), \ldots, \\
    &\phantom{\{} 5(y-0.4), \frac{25}{3}x(y-0.4), \ldots, \\
    &\phantom{\{} \frac{20}{7}(y-0.4)(y-0.6), \ldots \Bigg\}. \\
  \end{aligned}
  \end{equation}

\end{exmp}

\section{Algorithms and Timings}\label{sec:application}

At the beginning of this section is a restatement of the classical
BM algorithm with the notation established above.

\begin{alg}\label{BMalg}(BM Algorithm)

\textbf{Input}: A set of distinct points $\Xi=\{\xi^{(i)} : i=1,
\ldots, \mu\}\subset \mathbb{F}^d$ and
 a fixed monomial ordering $\prec$.\\
\indent \textbf{Output}: The 3-tuple $(G,N,Q)$, where $G$ is the
reduced Gr\"{o}bner basis for $\mathcal{I}(\Xi)$ w.r.t. $\prec$, $N$
is the Gr\"{o}bner \'{e}scalier of $\mathcal{I}(\Xi)$ (the DRIMB for
(\ref{e:inter}) also) w.r.t. $\prec$, and $Q$ is a DRINB
 w.r.t. $\prec$ for (\ref{e:inter}).
\\
\indent \textbf{BM1.} Start with lists $G=[\ ],N=[\ ],Q=[\ ],
L=[1]$, and a matrix
$B=(b_{ij})$ over $\mathbb{F}$ with $\mu$ columns and zero rows initially.\\
\indent \textbf{BM2.} If $L=[\ ]$, return $(G,N,Q)$ and stop.
Otherwise, choose the monomial
$t=\mbox{min}_\prec L$, and delete $t$ from $L$.\\
\indent \textbf{BM3.} Compute the evaluation vector
$(t(\xi^{(1)}),\ldots,t(\xi^{(\mu)}))$, and reduce it against the
rows of $B$ to obtain
$$
(v_1,\ldots,v_\mu)=(t(\xi^{(1)}),\ldots,t(\xi^{(\mu)}))-\sum_i
a_i(b_{i1},\ldots,b_{i\mu}), \quad a_i\in \mathbb{F}.
$$
\indent \textbf{BM4.}. If $(v_1,\ldots,v_\mu)=(0,\ldots,0)$, then
append the polynomial $t-\sum_i a_iq_i$ to the list $G$, where $q_i$
is the $i$th element of $Q$. Remove from $L$ all the multiples of
$t$.
Continue with \textbf{BM2}.\\
\indent \textbf{BM5.} Otherwise
$(v_1,\ldots,v_\mu)\neq(0,\ldots,0)$, add $(v_1,\ldots,v_\mu)$ as a
new row to $B$ and $t-\sum_i a_iq_i$ as a new element to $Q$. Append
the monomial $t$ to $N$, and add to $L$ those elements of
$\{x_1t,\ldots,x_dt\}$ that are neither multiples of an element of
$L$ nor of $\mathrm{LM}_\prec(G):=\{\mathrm{LM}_\prec(g): g\in G\}$.
Continue with \textbf{BM2}.
\end{alg}

Like TBM algorithm in \cite{DCZW2009}, we have the following QTBM
algorithm for quasi-tower sets.

\begin{alg}\label{XTBM}(QTBM algorithm)

\textbf{Input}: A quasi-$x$-tower(quasi-$y$-tower) set $\Xi\subset
\mathbb{F}^2$ of $\mu$ points and a fixed monomial ordering
$\prec_{\mathrm{grlex}}$($\prec_{\mathrm{grevlex}}$) or
$\prec_{\mathrm{lex}}$($\prec_{\mathrm{invlex}}$).

\textbf{Output}: The 3-tuple $(G,N,Q)$, where $G$ is the reduced
Gr\"{o}bner basis for $\mathcal{I}(\Xi)$, $N$ is the Gr\"{o}bner
\'{e}scalier of $\mathcal{I}(\Xi)$, and $Q$ is a DRINB for \eqref{e:inter}.\\
\indent \textbf{QTBM1.} Construct lower set $S_x(\Xi)$($S_y(\Xi)$).
\\
\indent \textbf{QTBM2.} Compute the sets $N$ and $Q$ according to
Theorem \ref{t:main}.
\\
\indent \textbf{QTBM3.} Compute the set $L:=\{x\cdot t : t\in
N\}\bigcup\{ y\cdot t : t\in N\}\setminus N$.
\\
\indent \textbf{QTBM4.} Construct $\mu\times \mu$ matrix $C$ whose
$(h, k)$ entry is $\phi_h^x(u_k^x)$($\phi_h^y(u_k^y)$) where
$\phi_h^x$($\phi_h^y$), $u_k^x$($u_k^y$) are $h$th and $k$th
elements of $Q=\{\phi_{ij}^x(\phi_{ij}^y): (i, j)\in
S_x(\Xi)(S_y(\Xi))\}$ and $\Xi=\{u_{mn}^x(u_{mn}^y): (m, n)\in
S_x(\Xi)(S_y(\Xi))\}$ w.r.t. the increasing
$\prec_{\mathrm{invlex}}$($\prec_{\mathrm{lex}}$) on $(i, j)$ and
$(m, n)$ respectively.
\\
\indent \textbf{QTBM5.} Goto \textbf{BM2} with $N,Q,L,C$ for the
reduced Gr\"{o}bner basis $G$.
\end{alg}

\begin{exmp}
  We continue with Example \ref{ex:qtt} and \ref{ex:qtt2}. Example \ref{ex:qtt2} had
  obtained $N$ and $Q$ in (\ref{e:exN}) and (\ref{e:exQ}) respectively. According to
  \textbf{QTBM3} and \textbf{QTBM4}, we get
  $$
  L=\{x^7, x^5y, x^6y, x^4y^2, x^3y^3, x^2y^4, y^5, xy^5\}
  $$
  and
  $$
  C=\left(
      \begin{array}{cccc}
        1 & 1 & 1 & \cdots \\
        0 & 1 & 3 & \cdots \\
        0 & 0 & 1 & \cdots \\
        \vdots & \vdots & \vdots & \ddots \\
      \end{array}
    \right).
  $$
  Finally, \textbf{QTBM5} yields
  \begin{align*}
    G=\{(y-0.4)(y-0.6)(y-1.1)(y-1.45)(y-1.8), \ldots\}.
  \end{align*}

\end{exmp}

At the end, we will show the timings for the computations of
BM-problems on quasi-tower sets in finite prime fields
$\mathbb{F}_q$ with size $q$ w.r.t. monomial orderings
$\prec_{\mathrm{grlex}}$ and $\prec_{\mathrm{invlex}}$ respectively.
QTBM and BM algorithms were implemented on Maple 12 installed on a
laptop with 768 Mb RAM and 1.5 GHz CPU.

For field $\mathbb{F}_{37}$ and $\prec_{\mathrm{grlex}}$,

\vskip 0.3cm

\begin{center}
\begin{tabular}{l|llll}\hline
  \backslashbox{Algorithms}{$\#\Xi$} & 300     & 500   & 800  &  1200\\ \hline\hline
  QTBM      & 2.563 s & 8.773 s & 26.068 s & 63.772 s\\
  BM     & 20.690 s & 80.375 s  & 288.364 s & 901.246 s\\
  \hline
\end{tabular}
\end{center}

\vskip 0.3cm

For field $\mathbb{F}_{43}$ and $\prec_{\mathrm{invlex}}$,

\vskip 0.3cm

\begin{center}
\begin{tabular}{l|llll}\hline
  \backslashbox{Algorithms}{$\#\Xi$}     & 500      & 800       & 1000   & 1200    \\ \hline\hline
  QTBM               & 8.171 s & 25.607 s & 48.170 s & 70.802 s \\
  BM              & 114.996 s  & 375.610 s  & 881.988 s  & 1347.287 \\
  \hline
\end{tabular}
\end{center}

\bibliographystyle{elsarticle-num}
\bibliography{refDWZL2010}

\end{document}